\documentclass[a4paper,10pt]{article}
\usepackage{amsfonts}
\usepackage{amstext}
\usepackage{amsthm}
\usepackage{amsmath}
\usepackage{amssymb}
\usepackage{latexsym}

\usepackage[latin1]{inputenc}
\newcommand{\bl}{\hfill\rule{2mm}{2mm}}
\newcommand{\R}{\Bbb{R}}

\newtheorem{teor}{Theorem}[section]

\newtheorem{lema}{Lemma}[section]

\newcommand{\n}{\noindent}

\begin{document}

\title{Optimal $L^p$-Riemannian Gagliardo-Nirenberg inequalities}
\author{\textbf{Jurandir Ceccon \footnote{\textit{E-mail addresses}:
ceccon@mat.ufpr.br (J. Ceccon)}}\\ {\small\it Departamento de
Matem\'{a}tica, Universidade Federal do Paran\'{a},}\\
{\small\it Caixa Postal 019081, 81531-990, Curitiba, PR, Brazil}\\
\textbf{Marcos Montenegro \footnote{\textit{E-mail addresses}:
montene@mat.ufmg.br (M. Montenegro)}}\\ {\small\it Departamento de
Matem\'{a}tica, Universidade Federal de Minas Gerais,}\\ {\small\it
Caixa Postal 702, 30123-970, Belo Horizonte, MG, Brazil}}

\date{}

\maketitle

\markboth{abstract}{abstract}
\addcontentsline{toc}{chapter}{abstract}

\hrule \vspace{0,2cm}

\n {\bf Abstract}

Let $(M,g)$ be a compact Riemannian manifold of dimension $n \geq 2$
and $1< p \leq 2$. In this work we prove the validity of the optimal
Gagliardo-Nirenberg inequality

\[
\left( \int_M |u|^r dv_g \right)^{\frac{p}{r \theta}} \leq \left(
A_{opt} \int_M |\nabla u|_g^p dv_g + B \int_M |u|^p dv_g \right)
\left( \int_M |u|^q dv_g \right)^{\frac{p(1 - \theta)}{\theta q}}
\]

\n for a family of parameters $r$, $q$ and $\theta$. Our proof relies strongly on a new distance lemma which holds for $1 < p \leq 2$. In particular, we obtain Riemannian versions of $L^p$-Euclidean Gagliardo-Nirenberg inequalities of \cite{DPDo} and extend the optimal $L^2$-Riemannian Gagliardo-Nirenberg inequality of \cite{Bro} in a unified framework.

\vspace{0,5cm} \hrule\vspace{0.2cm}

\section{Introduction and the main result}

Best constants and sharp first-order Sobolev inequalities on compact
Riemannian manifolds have been extensively studied in the last few
decades and surprising results have been obtained by showing the
influence of the geometry on such problems. Particularly, the
arising of concentration phenomena has motivated the development of
new methods in analysis, see \cite{H}, \cite{Au2} and \cite{DH} for
a complete survey. Our main interest in this work is the study of
optimal $L^p$-Riemannian Gagliardo-Nirenberg inequalities.
Such inequalities are connected for instance to a priori
bounds of solutions of some elliptic equations, see \cite{DelDo}.

Let $1 < p < n$ and $1 \leq q < r \leq p^*$, where $p^* =
\frac{np}{n-p}$ is the Sobolev critical exponent. Denote by
$D^{p,q}(\R^n)$ the completion of $C_0^\infty(\R^n)$ under the norm

\[
||u||_{D^{p,q}(\R^n)} = \left( \int_{\R^n} |\nabla u|^p dx
\right)^{\frac{1}{p}} + \left( \int_{\R^n} |u|^q dx
\right)^{\frac{1}{q}}\ .
\]

\n The $L^p$-Euclidean Gagliardo-Nirenberg inequality states that
for some constant $A$,

\begin{equation}\label{dgne}
\left( \int_{\R^n} |u|^r dx \right)^{\frac{p}{r \theta}} \leq A
\left( \int_{\R^n} |\nabla u|^p dx \right) \left( \int_{\R^n} |u|^q
dx \right)^{\frac{p(1 - \theta)}{\theta q}}
\end{equation}

\n for all $D^{p,q}(\R^n)$, where $\theta = \theta(p,q,r) =
\frac{np(r - q)}{r(q(p - n) + np)} \in (0,1]$. According to
\cite{Ba}, up to the constant $A$, these inequalities are all equivalent for $p$ fixed. Designate by $A(p,q,r)$ the best constant associated to (\ref{dgne}). Then, the optimal
$L^p$-Euclidean Gagliardo-Nirenberg inequality yields

\[
\left( \int_{\R^n} |u|^r dx \right)^{\frac{p}{r \theta}} \leq
A(p,q,r) \left( \int_{\R^n} |\nabla u|^p dx \right) \left(
\int_{\R^n} |u|^q dx \right)^{\frac{p(1 - \theta)}{\theta q}}\ .
\]

\n In \cite{DPDo}, Del Pino and Dolbeault studied this inequality
for $p < q < \frac{p(n - 1)}{n -
p}$ and $r = \frac{p(q - 1)}{p - 1}$. They proved that the
extremal functions associated to $A(p,q) = A(p,q,\frac{p(q - 1)}{p -
1})$ are given precisely by $u(x)=\alpha w(\beta(x-x_0))$, where
$\alpha\in\R$, $\beta\not=0$, $x_0\in\R^{n}$ and

\begin{equation}\label{ofgn}
w(x) = \left( 1 + \frac{q - p}{p - 1}|x|^{\frac{p}{p - 1}}
\right)^{- \frac{p - 1}{q - p}},
\end{equation}

\n and using this fact, they found

\[
A(p,q)= \left(\frac{q - p}{p \sqrt{\pi}} \right)^p \left(
\frac{pq}{n(q - p)}\right) \left(\frac{np - q(n -
p)}{pq}\right)^{\frac{p}{r \theta}} \left( \frac{\Gamma(\frac{q(p -
1)}{q - p}) \Gamma(\frac{n}{2} + 1)}{\Gamma(\frac{p - 1}{p} \frac{np
- q(n - p)}{q - p}) \Gamma(\frac{n(p - 1)}{p} + 1)}
\right)^{\frac{p}{n}}\ .
\]

Let $(M,g)$ be a compact Riemannian manifold of dimension $n
\geq 2$. An easy partition-of-unity argument and (\ref{dgne}) lead to the existence of constants $A,B \in \R$ such
that

\begin{gather}\label{AB1}
\left( \int_M |u|^r dv_g \right)^{\frac{p}{r \theta}} \leq \left( A
\int_M |\nabla_g u|^p dv_g + B \int_M |u|^p dv_g \right) \left(
\int_M |u|^q dv_g \right)^{\frac{p(1 - \theta)}{\theta q}}
\tag{$I_{p,q,r}(A,B)$}
\end{gather}

\n for all $u \in D^{p,q}(M)$, where $D^{p,q}(M)$ stands for the
completion of $C^\infty(M)$ with respect to the norm

\[
||u||_{D^{p,q}(M)} = \left( \int_M |\nabla_g u|^p dv_g
\right)^{\frac{1}{p}}+ \left( \int_M |u|^p dv_g
\right)^{\frac{1}{p}} + \left( \int_M |u|^q dv_g
\right)^{\frac{1}{q}}\ .
\]

\n In this case, we say simply that \ref{AB1} is valid. Such
inequality is known as $L^p$-Riemannian Gagliardo-Nirenberg
inequality. As pointed out by Brouttelande \cite{Bro}, \ref{AB1} are all equivalent, in the validity sense, for $p$ fixed.
This assertion follows from an adaption of the proof of Theorem 1.1
of \cite{DHV}.

The first best constant associated to \ref{AB1} is defined by

\[
A_{opt} = \inf \{ A \in \R:\; \mbox{there exists} \hspace{0,18cm} B
\in \R \hspace{0,18cm} \mbox{such that} \hspace{0,18 cm}
I_{p,q,r}(A,B) \hspace{0,18cm} \mbox{is valid}\}
\]

\n and the corresponding optimal inequality on $D^{p,q}(M)$ is $I_{p,q,r}(A_{opt},B)$, i.e

\[
\left( \int_M |u|^r dv_g \right)^{\frac{p}{r \theta}} \leq \left(
A_{opt} \int_M |\nabla_g u|^p dv_g + B \int_M |u|^p dv_g \right)
\left( \int_M |u|^q dv_g \right)^{\frac{p(1 - \theta)}{\theta q}} \
.
\]

\n Using again the idea of the proof of Theorem 1.1 of \cite{DHV},
one may easily check that $A_{opt} = A(p,q,r)$. In particular, the
explicit value of  $A_{opt}$ is known in the case $q > p$ and $r =
\frac{p(q -1 )}{p - 1}$. Some special cases of optimal
$L^p$-Riemannian Gagliardo-Nirenberg inequalities have been studied
in the last years. For instance, the optimal $L^p$-Riemannian
Sobolev inequality (i.e. $I_{p,q,r}(A_{opt},B)$ with $q =
\frac{(n-1)p}{n-p}$ and $r = p^*$) was proved to be valid for $1 < p
\leq 2$, independently, by Aubin and Li \cite{AuLi} and Druet
\cite{D2}, with best constant $A_{opt}$ found by Aubin \cite{Au1}
and Talenti \cite{Ta}. The optimal $L^2$-Riemannian Nash inequality
(i.e. $I_{p,q,r}(A_{opt},B)$ with $p=2$, $q = 1$ and $r = 2$) was
proved to be valid by Humbert \cite{Hu} with best constant $A_{opt}$
found by Carlen and Loss \cite{CaLo}. More generally, Brouttelande
proved in \cite{Bro} the validity of $I_{p,q,r}(A_{opt},B)$ for $p =
2$, $q < r$ and $1 \leq q \leq 2 \leq r < 2 + \frac{2}{n}q$. In
general, we cannot hope the validity of optimal $L^p$-Riemannian
Gagliardo-Nirenberg inequalities for $p > 2$. In fact, let $(M,g)$
be a $n$-dimensional compact Riemannian manifold with scalar
curvature $Scal_g$ positive in some point $x_0$ of $M$. Assume $1 < p < n$, $p < q < \frac{p(n - 1)}{n - p}$ and $r = \frac{p(q-1)}{p-1}$. For each $\varepsilon > 0$, define

\[
u_\varepsilon (\exp_{x_0}(x)) = \eta(x) w_{\varepsilon}(x),
\]

\n where $\eta \in C^{\infty}_0(B_\delta)$ is a cutoff function with
$\delta > 0$ smaller than the radius of injectivity of $M$, $w_{\varepsilon}(x) = \varepsilon^{-\frac{n}{p^*}}
w(\frac{x}{\varepsilon})$ and $w$ is the extremal function given in (\ref{ofgn}). Arguing in the same spirit of Druet
\cite{D1}, after several straightforward computations, we find for
$p > \max\{2, 2q/3\}$,

\[
\frac{\frac{\left( \int_M |u_\varepsilon|^r dv_g \right)^{\frac{p}{r \theta}}}{\left( \int_M |u_\varepsilon|^q dv_g \right)^{\frac{p(1 - \theta)}{\theta q}}} -
A_{opt} \int_M |\nabla_g u_\varepsilon|^p dv_g}{\int_M |u_\varepsilon|^p dv_g} =
\]

\[
\frac{\frac{1}{6n} Scal_g(x_0) \left(A_{opt} (
I_1^{\frac{p(1-\theta)}{\theta q}} I_2 + \frac{p(1-\theta)}{\theta
q} I_1 I_4 I_3^{\frac{p(1-\theta)}{\theta q} - 1}) -
\frac{p}{r \theta} I_5\right) \varepsilon^2 +
o(\varepsilon^2)}{o(\varepsilon^2)} \rightarrow + \infty
\]

\n as $\varepsilon \rightarrow 0$, where $I_k$ denotes the following integrals

\[
I_1 = \int_{\R^n} w^q\; dx, \ \ I_2 = \int_{\R^n} |\nabla w|^p |x|^2\; dx, \ \ I_3 = \int_{\R^n} w^q |x|^2\; dx,
\]

\[
I_4 = \int_{\R^n} |\nabla w|^p \; dx, \ \ I_5 = \int_{\R^n} w^r |x|^2\; dx\ .
\]

\n When $r = p^*$, this result recovers the non-validity of \cite{D1} since the condition above on $p$ yields $2 < p < \frac{n + 2}{3}$.

Our main result states the validity of $I_{p,q,r}(A_{opt},B)$ for $1 < p \leq 2$ and a family of parameters $r$ and $q$.

\begin{teor}\label{tgno1}
Let $(M,g)$ be a compact Riemannian manifold of dimension $n \geq 2$.
If either $1 < p \leq 2$, $p < r$ and $1 \leq q < r < p^*$ or $p = r > 1$, $q \geq 1$ and
$\frac{p^2}{2} \leq q < p$, then $I_{p,q,r}(A_{opt},B)$ is valid.
\end{teor}

Riemannian versions of $L^p$-Euclidean Gagliardo-Nirenberg inequalities of \cite{DPDo} in the case $1 < p \leq
2$ are in particular included in our theorem. Theorem \ref{tgno1} also
extends the optimal $L^2$-Riemmanian Gagliardo-Nirenberg inequality
of \cite{Bro} for $r > 2$. Our proof is inspired in the
works of Druet \cite{D2} and Brouttelande \cite{Bro}. New technical difficulties however
arise for other values of $p$. We remark that the proof given in \cite{Bro} for $p=2$ does not extend directly to $p \neq 2$. An illustrative example of this restriction is Lemma 2.3. In order to surround this obstacle, we use a version of the distance lemma (see the third step) which works well in the case $1 < p \leq 2$. Through these ideas it is also possible to simplify a little the proof of $L^p$-Riemannian Sobolev inequality given in \cite{D2}.

The optimal $L^p$-Riemannian logarithmic Sobolev inequality states

\[
\int_M |u|^p \log(|u|^p) dv_g \leq \frac{n}{p} \log \left( {\cal A}_{opt} \int_M |\nabla_g u|^p dv_g + B \right)
\]

\n for all $u \in H^{1,p}(M)$ such that $\|u\|_{L^p} = 1$. The validity of this inequality is an open question for $1 < p < 2$. We remark that Theorem \ref{tgno1} provides a tool to investigate the validity of this inequality. In fact, consider the second best constant associated to \ref{AB1} given by

\[
B(p,q,r) = \inf \{ B \in \R:\; I_{p,q,r}(A_{opt},B) \hspace{0,18cm} \mbox{is valid}\}\ .
\]

\n Let $q > p$ and $r = \frac{p(q-1)}{p-1}$. In this case, the value of $A_{opt}$ is $A(p,q)$, see \cite{DelDo}. If $B(p,q,r)$ is bounded on $q$, then, taking $q \downarrow p$ in $I_{p,q,r}(A_{opt},B(p,q,r))$, one obtains the desired logarithmic inequality as a limit case of optimal $L^p$-Riemannian Gagliardo-Nirenberg inequalities. Note, however, that the a priori estimates problem for the second best constant is in general quite delicate.

\section{Proof of Theorem \ref{tgno1}}

We proceed by contradiction as in
\cite{D2} and \cite{Bro}. So, for each $\alpha > 0$, one has

\begin{equation}\label{3dha}
\nu_\alpha = \inf_{u \in E} J_\alpha (u) < A(p,q,r)^{-1},
\end{equation}

\n where $E = \{ u \in D^{p,q}(M):\; ||u||_{L^r(M)} = 1 \}$
and

\[
J_\alpha(u) = \left( \int_M |\nabla_g u|^p dv_g + \alpha \int_M
|u|^p dv_g \right) \left( \int_M |u|^q dv_g \right)^{\frac{p(1 -
\theta)}{\theta q}}\ .
\]

\n Since $J_\alpha$ is of class $C^1$, using standard variational
arguments, we find a minimizer $u_\alpha \in E$ of
$J_\alpha$, i.e.

\begin{equation}\label{3iha}
J_\alpha(u_\alpha) = \nu_\alpha = \inf_{u \in E} J_\alpha(u)\
.
\end{equation}

\n One may assume $u_\alpha \geq 0$ since $\nabla_g |u_\alpha| = \pm
\nabla_g u_\alpha$. Clearly, $u_\alpha$ satisfies the Euler-Lagrange
equation

\begin{equation} \label{3ep}
A_\alpha \Delta_{p,g} u_\alpha + \alpha A_\alpha u_\alpha^{p - 1} +
\frac{1 - \theta}{\theta} B_\alpha u_\alpha^{q - 1} = \mu_\alpha
u_\alpha^{r - 1},
\end{equation}

\n where $\Delta_{p,g} = -{\rm div}_g(|\nabla_g|^{p-2} \nabla_g)$ is
the $p$-Laplace operator on $M$,

\[
A_\alpha = \left(\int_M u_\alpha^q dv_g \right)^{\frac{p( 1 -
\theta)}{\theta q}},
\]

\[
B_\alpha = \left(\int_M |\nabla_g u_\alpha|^p dv_g + \alpha \int_M
u_\alpha ^p dv_g \right) \left(\int_M u_\alpha ^q dv_g
\right)^{\frac{p(1 - \theta)}{\theta q} - 1}
\]

\n and

\[
\mu_\alpha = \frac{1}{\theta}\nu_\alpha\ .
\]

\n By a regularity result due to Tolksdorf \cite{To}, it follows
that $u_\alpha$ is of class $C^1$.

In the following, we divide the proof into three steps. Several possibly
different positive constants independent of $\alpha$ are denoted by
$c$ and $c_i$. We also assume, without loss of generality, that the
radius of injectivity of $M$ is greater than $2$. \\

\n {\bf First step:} The following convergences\\

\n (a) $A_\alpha\int_M |\nabla_g u_\alpha|^p dv_g \rightarrow
A(p,q,r)^{-1}$,\\

\n (b) $\mu_\alpha \rightarrow \frac{A(p,q,r)^{-1}}{\theta}$,\\

\n (c) $\alpha A_\alpha \int_M u_\alpha^p dv_g \rightarrow 0$,\\

\n (d) $B_\alpha \int_M u_\alpha^q dv_g \rightarrow A(p,q,r)^{-1}$\\

\n hold as $\alpha \rightarrow \infty$.\\

\n {\bf Proof of (a):}  It follows from (\ref{3dha}) and
(\ref{3iha}) that

\[
\alpha A_\alpha \int_M u_\alpha^p dv_g < A(p,q,r)^{-1},
\]

\n so that

\begin{equation}\label{3dtf}
A_\alpha \int_M u_\alpha^p dv_g \rightarrow 0\ .
\end{equation}

\n From $A_{opt} = A(p,q,r)$ and (\ref{3dha}), one arrives at

\[
1 \leq A_\alpha \left( (A(p,q,r) + \varepsilon) \int_M |\nabla_g
u_\alpha|^p dv_g + B_\varepsilon \int_M u_\alpha^p dv_g \right)
\]

\n for some constant $B_\varepsilon$ independent of $\alpha$.
Letting $\alpha \rightarrow \infty$ and $\varepsilon \rightarrow 0$
and using (\ref{3dtf}), one has

\[
\liminf_{\alpha \rightarrow \infty} \left( A_\alpha \int_M |\nabla_g
u_\alpha|^p dv_g \right) \geq A(p,q,r)^{-1}\ .
\]

\n Noting that $A_\alpha\int_M |\nabla_g u_\alpha|^p dv_g < A(p,q,r)^{-1}$
for all $\alpha$, we conclude the part (a).\bl \\

\n {\bf Proof of (b):} This assertion follows directly from

\[
A_\alpha \int_M |\nabla_g u_\alpha|^p dv_g \leq \nu_\alpha \leq
A(p,q,r)^{-1}
\]

\n and the part (a).\bl \\

\n {\bf Proof of (c):} Taking the limit in

\[
A_\alpha \int_M |\nabla_g u_\alpha|^p dv_g + \alpha A_\alpha \int_M
u_\alpha^p dv_g < A(p,q,r)^{-1}
\]

\n and again using (a), we find (c).\bl \\

\n {\bf Proof of (d):} From definitions of $A_\alpha$ and
$B_\alpha$, one has

\[
B_\alpha \int_M u_\alpha^q dv_g = A_\alpha \int_M |\nabla_g
u_\alpha|^p dv_g + \alpha A_\alpha \int_M u_\alpha^p dv_g\ .
\]

\n So, the part (d) follows directly from (a) and (c).\bl\\

Let $x_\alpha \in M$ be a maximum point of $u_\alpha$, i.e

\begin{equation}\label{3ix}
u_\alpha(x_\alpha) = ||u_\alpha||_{L^\infty(M)}\ .
\end{equation}

\n {\bf Second step:} Let $(c_\alpha)_{\alpha > 0}$ be a family of
positive numbers such that $\frac{A_\alpha^{\frac{1}{p}}
||u_\alpha||_{L^\infty(M)}^{1 - \frac{r}{p}} }{c_\alpha} \rightarrow
0$. Then

\[
\int_{B(x_\alpha,c_\alpha)} u_\alpha^r dv_g \rightarrow 1 \ .
\]

The proof of this step is fairly technical and long. In order to make it clear we state two lemmas.

\begin{lema} \label{3lmz}
Let $D_\alpha = B(x_\alpha, A_\alpha^{\frac{1}{p}}
||u_\alpha||_{L^\infty(M)}^{1 - \frac{r}{p}})$. There exists a
constant $c > 0$ such that

\[
\frac{\int_{D_\alpha} u_\alpha^q dv_g}{\int_M u_\alpha ^q dv_g} \geq
c
\]

\n for $\alpha$ large.
\end{lema}

\n {\bf Proof of Lemma \ref{3lmz}.} By (\ref{3iha}) and (b) of the
first step, one has

\[
\alpha A_\alpha \int_M u_\alpha^p dv_g \leq c \int_M u_\alpha^r dv_g
\leq c ||u_\alpha||_{L^\infty(M)}^{r - p} \int_M u_\alpha^p dv_g,
\]

\n so that

\begin{equation}\label{fim31}
\alpha A_\alpha ||u_\alpha||_{L^\infty(M)}^{p - r} \leq c \ .
\end{equation}

\n In particular, $A_\alpha ||u_\alpha||_{L^\infty(M)}^{p - r}
\rightarrow 0$. For $ x \in B(0,1)$, define

\[
\begin{array}{l}
h_\alpha(x) = g(\exp_{x_\alpha} (A_\alpha^{\frac{1}{p}}
||u_\alpha||_{L^\infty(M)}^{1 - \frac{r}{p}}\; x)), \vspace{0,3cm}\\
\varphi_\alpha(x) = ||u_\alpha||_{L^\infty(M)}^{-1}
u_\alpha(\exp_{x_\alpha}(A_\alpha^{\frac{1}{p}}
||u_\alpha||_{L^\infty(M)}^{1 - \frac{r}{p}}\; x))\ .
\end{array}
\]

\n Clearly, $\varphi_\alpha$ satisfies

\[
\Delta_{p,h_\alpha} \varphi_\alpha + \alpha A_\alpha
||u_\alpha||_{L^\infty(M)}^{p - r} \varphi_\alpha^{p - 1} + \frac{1
- \theta}{\theta} B_\alpha ||u_\alpha||_{L^\infty(M)}^{q - r}
\varphi_\alpha^{q - 1} = \mu_\alpha \varphi_\alpha^{r - 1}\ .
\]

\n Arguing as in (\ref{fim31}), we find

\[
B_\alpha ||u_\alpha||_{L^\infty(M)}^{q - r} \leq c,
\]

\n so that all coefficients of the equation above are bounded. So,
from $||\varphi_\alpha||_{L^\infty(M)} \leq 1$ and an a priori estimate due to Tolksdorf \cite{To}, one concludes that $\varphi_\alpha$ converges to $\varphi$ in $C^1_{loc}(B(0,1)))$. In particular,
$\varphi \not\equiv 0$ since $\varphi_\alpha(0) = 1$ for all
$\alpha$. We claim that

\begin{equation}\label{3dqp}
1 \leq ||u_\alpha||_{L^\infty(M)}^{r - q} A_\alpha^{\frac{\theta
q}{p(1 - \theta)}} \leq c
\end{equation}

\n for some $c > 0$ independent of $\alpha$. In fact, by the Cartan
expansion of $h_\alpha$, we have

\[
\int_{B(0,1)} \varphi_\alpha^r dv_{h_\alpha} \geq c
\int_{B(0,\delta)} \varphi_\alpha^r dx \rightarrow c
\int_{B(0,\delta)} \varphi^r dx > 0
\]

\n for each $0 < \delta < 1$. Thus, from

\[
\int_{B(0,1)} \varphi_\alpha^r dv_{h_\alpha} =
\left(||u_\alpha||_{L^\infty(M)}^{\frac{np + rp -nr}{p}}
A_\alpha^{\frac{n}{p}} \right)^{-1} \int_{D_\alpha} u_\alpha^r dv_g,
\]

\begin{equation}\label{dt}
\frac{\theta q}{p(1 - \theta)} = \frac{n(r - q)}{np + rp - nr}
\end{equation}

\n and $||u_\alpha||_{L^r(M)} = 1$, we find

\[
||u_\alpha||_{L^\infty(M)}^{r - q} A_\alpha^{\frac{\theta q}{p(1 -
\theta)}} = ||u_\alpha||_{L^\infty(M)}^{\frac{np + rp -
nr}{n}\frac{\theta q}{p(1 - \theta)}} A_\alpha^{ \frac{\theta q}{p(1
- \theta)}} = \left(||u_\alpha||_{L^\infty(M)}^{\frac{np + rp -
nr}{p}} A_\alpha ^{\frac{n}{p}} \right)^{\frac{\theta q}{n(1 -
\theta)}} \leq c\ .
\]

\n On the other hand, from the definition of $A_\alpha$, we have

\[
1 = \int_M u_\alpha^r dv_g \leq ||u_\alpha||_{L^\infty(M)}^{r - q}
A_\alpha^{\frac{\theta q}{p(1- \theta)}},
\]

\n so that (\ref{3dqp}) holds. By (\ref{dt}), we may write

\[
\int_{B(0,1)} \varphi_\alpha^q dv_{h_\alpha} =
||u_\alpha||_{L^\infty(M)}^{\frac{nr -pq -np}{p}}
A_\alpha^{-\frac{n}{p} + \frac{\theta q}{p( 1 - \theta)}}
\frac{\int_{D_\alpha} u_\alpha^q dv_g}{\int_M u_\alpha ^q dv_g} =
\left(||u_\alpha||_{L^\infty(M)}^{r - q} A_\alpha ^{\frac{\theta
q}{p(1 - \theta)}}\right)^{\frac{nr -pq - np}{p(r - q)}}
\frac{\int_{D_\alpha} u_\alpha^q dv_g}{\int_M u_\alpha^q dv_g}\ .
\]

\n The proof then follows from (\ref{3dqp}) and $\varphi_\alpha
\rightarrow \varphi$ in $C_{loc}(B(0,1))$.\bl\\

\begin{lema}\label{3lmi}
Let $(c_\alpha)_{\alpha > 0}$ be a family of positive numbers such
that $\frac{A_\alpha^{\frac{1}{p}} ||u_\alpha||_{L^\infty(M)}^{1 -
\frac{r}{p}} }{c_\alpha} \rightarrow 0$. Then

\[
\frac{\int_{B(x_\alpha,c_\alpha)} u_\alpha^q dv_g}{\int_M u_\alpha^q
dv_g} \rightarrow 1\ .
\]
\end{lema}

\n {\bf Proof of Lemma \ref{3lmi}.} Let $\eta \in C_0^1(\R)$ be a
cutoff function such that $\eta = 1$ on $[0,\frac{1}{2}]$, $\eta =
0$ on $[1,\infty)$ and $0 \leq \eta \leq 1$. Define
$\eta_{\alpha,k}(x) = \eta(c_\alpha^{-1} d_g(x,x_\alpha))^{\tau^k}$,
where $\tau = \frac{p^*}{q}$. Taking $u_\alpha \eta_{\alpha,k}^r$ as
a test function in (\ref{3ep}), one has

\[
A_\alpha \int_M |\nabla_g u_\alpha|^p \eta_{\alpha,k}^r dv_g +
A_\alpha \int_M |\nabla_g u_\alpha|^{p - 2} u_\alpha \nabla_g
u_\alpha \cdot \nabla_g(\eta_{\alpha,k}^r) dv_g + \alpha A_\alpha
\int_M u_\alpha^p \eta_{\alpha, k}^r dv_g
\]

\begin{equation}\label{3sc3}
+ \frac{1 - \theta}{\theta} B_\alpha \int_M u_\alpha^q
\eta_{\alpha,k}^r dv_g = \mu_\alpha \int_M u_\alpha^r
\eta_{\alpha,k}^r dv_g\ .
\end{equation}

\n By $|\nabla_g \eta_{\alpha,k}| \leq \frac{c}{c_\alpha}$ and
(\ref{3dqp}),

\[
A_\alpha \int_{B(x_\alpha, c_\alpha)} u_\alpha^p |\nabla_g
\eta_{\alpha,k}|^p dv_g \leq c \frac{A_\alpha}{c_\alpha^p}
\int_{B(x_\alpha,c_\alpha)} u_\alpha^p dv_g \leq c
\frac{A_\alpha}{c_\alpha^p} \left(\int_M u_\alpha^r dv_g
\right)^{\frac{p}{r}} \left( \int_{B(x_\alpha, c_\alpha)} dv_g
\right)^{1 - \frac{p}{r}}
\]

\begin{equation}\label{3dqs3}
\leq c \frac{||u_\alpha||_{L^\infty(M)}^{( q - r)\frac{p(1 -
\theta)}{\theta q}}}{c_\alpha^{\frac{rp - nr + np}{r}}} = c \left(
\frac{||u_\alpha||_{L^\infty(M)}^{- \frac{r}{n}}}{c_\alpha}
\right)^{\frac{rp - nr + np}{r}} \leq c \left(
\frac{A_\alpha^{\frac{1}{p}} ||u_\alpha||_{L^\infty(M)}^{1 -
\frac{r}{p}}}{c_\alpha} \right)^{\frac{rp - nr + np}{r}}\ .
\end{equation}

\n So, by (a) of the first step, $r < p^*$ and the assumption of the
lemma, one concludes that

\[
A_\alpha \int_M |\nabla_g u_\alpha|^{p - 1} u_\alpha |\nabla_g
\eta_{\alpha,k}| dv_g \rightarrow 0\ .
\]

\n Using (b), (c) and (d) of the first step in
(\ref{3sc3}), we find

\[
\lim_{\alpha \rightarrow \infty}\left(A_\alpha \int_M |\nabla_g
u_\alpha|^p \eta_{\alpha, k}^r dv_g \right) + \frac{1 -
\theta}{\theta} A(p,q,r)^{-1} \lim_{\alpha \rightarrow \infty}
\left(\frac{\int_M u_\alpha^q \eta_{\alpha,k}^r dv_g}{\int_M
u_\alpha^q dv_g} \right)
\]

\begin{equation}\label{3il}
=\frac{A(p,q,r)^{-1}}{\theta} \lim_{\alpha \rightarrow \infty}
\int_M u_\alpha^r \eta_{\alpha, k}^r dv_g\ .
\end{equation}

\n On the other hand, for each $\varepsilon > 0$ there exists
$B_\varepsilon > 0$ such that

\[
\left( \int_M u_\alpha^r \eta^r_{\alpha,k} dv_g \right)^{\frac{p}{r
\theta}} \leq \left( (A(p,q,r) + \varepsilon) \int_M
|\nabla_g(u_\alpha \eta_{\alpha,k})|^p dv_g + B_\varepsilon \int_M
u_\alpha^p \eta_{\alpha,k}^p dv_g \right) \left(\int_M u_\alpha^q
\eta^q_{\alpha,k} dv_g \right)^{\frac{p(1 - \theta)}{\theta q}}\ .
\]

\n So, by the definition of $A_\alpha$, one has

\[
\left( \int_M u_\alpha^r \eta^r_{\alpha,k} dv_g \right)^{\frac{p}{r
\theta}} \leq (A(p,q,r) + \varepsilon) \left( \int_M |\nabla_g
u_\alpha|^p \eta_{\alpha,k}^p dv_g\right) \left(\int_M u_\alpha^q
\eta_{\alpha,k}^q dv_g \right)^{\frac{p(1 - \theta)}{\theta q}}
\]

\begin{equation}\label{3dqs5}
+ c A_\alpha \int_M |\nabla_g u_\alpha|^{p - 1} \eta_{\alpha,k}^{p -
1} u_\alpha |\nabla_g \eta_{\alpha,k}| dv_g + c A_\alpha \int_M
u_\alpha^p |\nabla_g \eta_{\alpha,k}|^p dv_g + c A_\alpha \int_M
u_\alpha^p dv_g\ .
\end{equation}

\n Therefore, by (\ref{3dqs3}) and (c) of the first step,

\begin{equation}\label{3il2}
\lim_{\alpha \rightarrow \infty}\left( \int_M u_\alpha^r
\eta_{\alpha,k}^r dv_g \right)^{\frac{p}{r \theta}} \leq A(p,q,r)
\lim_{\alpha \rightarrow \infty}\left(A_\alpha \int_M |\nabla_g
u_\alpha|^p \eta_{\alpha,k}^p dv_g \right) \lim_{\alpha \rightarrow
\infty} \left(\frac{\left(\int_M u_\alpha^q \eta_{\alpha,k}^q dv_g
\right)^{\frac{p(1 - \theta)}{\theta q}}}{A_\alpha} \right)\ .
\end{equation}

\n Let

\[
X_k = A(p,q,r) \lim_{\alpha \rightarrow \infty} \left( A_\alpha
\int_M |\nabla_g u_\alpha|^p \eta_{\alpha,k}^r dv_g \right), \ \ Y_k
= A(p,q,r) \lim_{\alpha \rightarrow \infty} \left( A_\alpha \int_M
|\nabla_g u_\alpha|^p \eta_{\alpha,k}^p dv_g \right)
\]

\[
Z_k = \lim_{\alpha \rightarrow \infty} \int_M \eta_{\alpha,k}^r
u_\alpha^r dv_g,\ \ \lambda_k = \lim_{\alpha \rightarrow
\infty}\frac{\int_M u_\alpha^q \eta_{\alpha,k}^r dv_g}{\int_M
u_\alpha^q dv_g},\ \ \tilde{\lambda}_k = \lim_{\alpha \rightarrow
\infty} \frac{\int_M ( \eta_{\alpha,k} u_\alpha)^q dv_g}{\int_M
u_\alpha^q dv_g}\ .
\]

\n Then, (\ref{3il}) and (\ref{3il2}) may be written as

\begin{equation}\label{3sc4}
\lambda_k = \frac{1}{1 - \theta}(Z_k - \theta X_k)
\end{equation}

\n and

\[
Z_k^{\frac{p}{r \theta}} \leq Y_k \tilde{\lambda}_k^{\frac{p(1 -
\theta)}{\theta q}},
\]

\n or equivalently,

\[
\lambda_k = \frac{1}{ 1 - \theta} Y_k^{\frac{r \theta}{p(1 -
\theta)}} \left( Z_k Y_k^{-\frac{r \theta}{p(1 - \theta)}}
\tilde{\lambda}_k^{-\frac{r}{q}} - \theta X_k Y_k^{-\frac{r
\theta}{p(1 - \theta)}} \tilde{\lambda}_k^{-\frac{r}{q}} \right)
\tilde{\lambda}_k^{\frac{r}{q}}
\]

\n and

\[
Y_k^{-\frac{r\theta}{p(1 - \theta)}}
\tilde{\lambda}_k^{-\frac{r}{q}} \leq Z_k^{- \frac{1}{ 1 - \theta}}\
.
\]

\n These relations lead us to

\begin{equation}\label{3dp2}
\lambda_k \leq \frac{1}{1 - \theta} Y_k^{\frac{r \theta}{p(1 -
\theta)}} \left( Z_k^{1 - \frac{1}{1 - \theta}} - \theta X_k Z_k^{-
\frac{1}{1 - \theta}} \right) \tilde{\lambda}_k^{\frac{r}{q}}.
\end{equation}

\n Define $f(x,z) = z^{1 - \frac{1}{1 - \theta}} - \theta x z^{-
\frac{1}{1 - \theta}}$. It follows easily that $f$ is non-increasing
in $z$ for $0 < x \leq z$ and non-decreasing in $z$ for $x \geq z >
0$. Noting that $Z_k = \theta X_k + (1 - \theta) \lambda_k$ implies
that either $\lambda_k \leq Z_k \leq X_k$ or $X_k \leq Z_k \leq
\lambda_k$, then $f(X_k,Z_k) \leq f(X_k,X_k)$. So, by
(\ref{3dp2}),

\[
\lambda_k \leq \left(Y_k^{\frac{r}{p}} X_k^{-1}
\right)^{\frac{\theta}{1 - \theta}} \tilde{\lambda}_k^{\frac{r}{q}}\
.
\]

\n On the other hand, from (a) of the first step, one has

\[
Y_k^{\frac{r}{p}} \leq \left[A(p,q,r)^{\frac{p}{\theta}}
\lim_{\alpha \rightarrow \infty} \left( A_\alpha \int_M
\eta_{\alpha,k}^r |\nabla u_\alpha|_g^p dv_g \right)^{\frac{p}{r}}
\lim_{\alpha \rightarrow \infty} \left( A_\alpha \int_M |\nabla
u_\alpha|_g^p dv_g\right)^{\frac{r - p}{r}} \right]^{\frac{r}{p}} =
X_k\ .
\]

\n So,

\begin{equation}\label{3dl2}
\lambda_k \leq \tilde{\lambda}_k^{\frac{r}{q}}\ .
\end{equation}

\n By Lemma \ref{3lmz}, for any $k$, we have

\[
0 < c < \lim_{\alpha \rightarrow \infty} \frac{\int_{B(x_\alpha,
\frac{c_\alpha}{2})} u_\alpha^q dv_g}{\int_M u_\alpha^q dv_g} \leq
\lim_{\alpha \rightarrow \infty} \frac{\int_M u_\alpha^q
\eta_{\alpha,k}^r dv_g}{\int_M u_\alpha^q dv_g} = \lambda_k\ .
\]

\n Since $q < r$ and $\tau = \frac{p^*}{q}$, one concludes easily
from (\ref{3dl2}) that

\[
\lambda_k \leq \tilde{\lambda}_k \leq \lambda_{k - 1}\ .
\]

\n These inequalities lead us to

\[
0 < c < \cdots \leq \lambda_{k + 1} \leq \tilde{\lambda}_{k + 1}
\leq \lambda_k \leq \cdots \leq \lim_{\alpha \rightarrow \infty}
\frac{\int_{B(x_\alpha, c_\alpha)}u_\alpha^q dv_g}{\int_M u_\alpha^q
dv_g}\ .
\]

\n Noting that

\[
\lambda_0 = \lim_{\alpha \rightarrow \infty}
\frac{\int_{B(x_\alpha,c_\alpha)} \eta_{\alpha,0}^r u_\alpha^q
dv_g}{\int_M u_\alpha^q dv_g} \leq \lim_{\alpha \rightarrow \infty}
\frac{\int_{B(x_\alpha,c_\alpha)} u_\alpha^q dv_g}{\int_M u_\alpha^q
dv_g} \leq 1,
\]

\n it follows from (\ref{3dl2}) that $c <
(\lambda_0^{\frac{r}{q}})^k \leq 1$ for all $k$. Therefore, $\lambda_0 =
1$ and this concludes the
proof.\bl \\

\n {\bf Proof of the second step:} From (\ref{3dqp}) and Lemma
\ref{3lmi}, we find

\[
\int_{M \backslash B(x_\alpha,c_\alpha)} u_\alpha^r dv_g \leq
A_\alpha^{\frac{\theta q}{p(1 - \theta)}}
||u_\alpha||_{L^\infty(M)}^{r - q} \frac{\int_{M \backslash
B(x_\alpha,c_\alpha)} u_\alpha^q dv_g}{\int_M u_\alpha^q dv_g} \leq
c \frac{\int_{M \backslash B(x_\alpha,c_\alpha)} u_\alpha^q
dv_g}{\int_M u_\alpha^q dv_g} \rightarrow 0\ .
\]

\n So,

\[
\int_{B(x_\alpha,c_\alpha)} u_\alpha^r dv_g \rightarrow 1
\]

\n since $||u_\alpha||_r = 1$.\bl \\

The next step states a priori estimates of $u_\alpha$ in two
distinct cases: $p < r$ and $p = r$. Note that here it arises a
discontinuity interesting phenomenon with respect to $p$. Such estimates involve the
weight $A_\alpha$ and this dependence is essential in proof of Theorem \ref{tgno1}.\\

\n {\bf Third step:} There exists a constant $c>0$ such that\\

\n (a)\ $A_\alpha^{-\frac{1}{p}} d_g(x,x_\alpha)
u_\alpha(x)^{\frac{r -
p}{p}} \leq c$ if $p < r$,\\

\n (b)\ $ d_g(x,x_\alpha) u_\alpha(x)^{\frac{p}{n}} \leq c$ if $p =
r$

\vspace{0,3cm}

\n for all $x \in B(x_\alpha,1)$, where $d_g$ stands for the
distance
with respect to $g$.\\

\n {\bf Proof of the third step:} Suppose, by contradiction, that
the assertion (a) is false.  Set $w_\alpha(x) =
A_\alpha^{-\frac{1}{r - p}} d_g(x,x_\alpha)^{\frac{p}{r-p}}
u_\alpha(x)$. Let $y_\alpha \in \overline{B(x_\alpha,1)}$ be a
maximum point of $w_\alpha$. Then, $w_\alpha(y_\alpha) =
||w_\alpha||_{L^\infty(B(x_\alpha,1))} \rightarrow \infty$. We first
show that

\begin{equation}\label{3int}
B(y_\alpha,A_\alpha^{\frac{1}{p}}
u_\alpha(y_\alpha)^{\frac{p-r}{p}}) \cap B(x_\alpha,
w_\alpha(y_\alpha)^\nu A_\alpha^{\frac{1}{p}}
u_\alpha(y_\alpha)^{\frac{p-r}{p}}) = \emptyset
\end{equation}

\n for some $\nu > 0$ and $\alpha$ large. Note that this fact is
clearly implied by

\[
d_g(x_\alpha,y_\alpha) \geq A_\alpha^{\frac{1}{p}}
u_\alpha(y_\alpha)^{\frac{p - r}{p}} + w_\alpha(y_\alpha)^\nu
A_\alpha ^{\frac{1}{p}} ||u_\alpha||_{L^\infty(M)}^{\frac{p - r}{p}}
\ .
\]

\n But this inequality is equivalent to

\[
w_\alpha(y_\alpha)^{\frac{r - p}{p} - \nu} \geq
w_\alpha(y_\alpha)^{- \nu} + u_\alpha(y_\alpha)^{\frac{r - p}{p}}
||u_\alpha||_{L^\infty(M)}^{\frac{p - r}{p}}\ .
\]

\n Taking $0 < \nu < \frac{r - p}{p}$ we see that the inequality
above holds for $\alpha$ large since $w_\alpha(y_\alpha)^{-\nu}
\rightarrow 0$ and $w_\alpha(y_\alpha)^{\frac{r - p}{p} - \nu}
\rightarrow \infty$. This ends the proof of (\ref{3int}). Note that
$A_\alpha^{\frac{1}{p}} u_\alpha(y_\alpha)^{\frac{p-r}{p}}
\rightarrow 0$ since $||w_\alpha||_{L^\infty(B(x_\alpha,1))} \rightarrow
\infty$. For $x \in B(0,1)$, define

\[
\begin{array}{l}
h_\alpha(x) = g(\exp_{y_\alpha}(A_\alpha^{\frac{1}{p}}
u_\alpha(y_\alpha)^{\frac{p-r}{p}} x)) \vspace{0,2cm}\\
\psi_\alpha(x) = u_\alpha(y_\alpha)^{-1}
u_\alpha(\exp_{y_\alpha}(A_\alpha^{\frac{1}{p}}
u_\alpha(y_\alpha)^{\frac{p-r}{p}} x))\ .
\end{array}
\]

\n Clearly, $\psi_\alpha$ satisfies

\begin{equation}\label{3sc5}
\Delta_{p,h_\alpha} \psi_\alpha + \alpha A_\alpha
u_\alpha(y_\alpha)^{p - r} \psi_\alpha^{p - 1} + \frac{1 -
\theta}{\theta} B_\alpha u_\alpha(y_\alpha)^{q - r} \psi_\alpha^{q -
1} = \mu_\alpha \psi_\alpha^{r - 1}\ .
\end{equation}

\n Let $0 < c < 1$. We claim that $u_\alpha(y_\alpha) \geq c
u_\alpha(x)$ for all $x \in B(y_\alpha, A_\alpha^{\frac{1}{p}}
u_\alpha(y_\alpha)^{\frac{p-r}{p}})$ and $\alpha$ large. In fact,
set $c_0 = c^{\frac{r-p}{p}}$. Then, $A_\alpha^{\frac{1}{p}}
u_\alpha(y_\alpha)^{\frac{p - r}{p}} \leq (1-c_0)
d_g(x_\alpha,y_\alpha)$ for $\alpha$ large, so that

\[
d_g(x_\alpha,x) \geq d_g(x_\alpha,y_\alpha) - d_g(x,y_\alpha) \geq
d_g(x_\alpha,y_\alpha) - A_\alpha^{\frac{1}{p}}
u_\alpha(y_\alpha)^{\frac{p - r}{p}} \geq c_0
d_g(x_\alpha,y_\alpha)\ .
\]

\n So, the desired estimate follows directly from

\[
A_\alpha^{-\frac{1}{r - p}} d_g(x_\alpha,y_\alpha)^{\frac{p}{r - p}}
u_\alpha(y_\alpha) = w_\alpha(y_\alpha) \geq w_\alpha(x) =
A_\alpha^{-\frac{1}{r - p}} d_g(x_\alpha,x)^{\frac{p}{r - p}}
u_\alpha(x) \ .
\]

\n In particular,

\begin{equation}\label{3lu}
||\psi_\alpha||_{L^\infty(B(0,1))} \leq c
\end{equation}

\n for $\alpha$ large. By (\ref{3sc5}),

\[
\int_{B(0,1)} |\nabla_{h_\alpha} \psi_\alpha|^{p - 2}
\nabla_{h_\alpha} \psi_\alpha \cdot  \nabla_{h_\alpha} \phi\;
dv_{h_\alpha} \leq c \int_{B(0,1)} (\psi_\alpha^{r - p})
\psi_\alpha^{p - 1} \phi\; dv_{h_\alpha}
\]

\n for all positive test function $\phi \in C_0^1(B(0,1))$. So, by
Moser's iterative scheme, (\ref{3dqp}) and (\ref{dt}), it follows
that

\[
1 = \sup_{B(0,\frac{1}{2})} \psi_\alpha \leq c \int_{B(0,1)}
\psi_\alpha^r dv_{h_\alpha} = c\left( A_\alpha^{\frac{\theta q}{p(1
- \theta)}} u_\alpha(y_\alpha)^{r - q} \right)^{-\frac{n(1 -
\theta)}{\theta q}} \int_{\tilde{D}_\alpha} u_\alpha^r dv_g
\]

\[
 \leq c
\left(\frac{||u_\alpha||_{L^\infty(M)}}{u_\alpha(y_\alpha)}
\right)^{\frac{np - rn + pr}{p}} \int_{\tilde{D}_\alpha} u_\alpha^r
dv_g,
\]

\n where $\tilde{D}_\alpha = B(y_\alpha, A_\alpha^{\frac{1}{p}}
u_\alpha(y_\alpha)^{\frac{p-r}{p}})$. This last inequality may be
rewritten as

\begin{equation}\label{3dcm}
0 < c \leq m_\alpha^\sigma \int_{\tilde{D}_\alpha} u_\alpha^r dv_g,
\end{equation}

\n where $m_\alpha =
\frac{||u_\alpha||_{L^\infty(M)}}{u_\alpha(y_\alpha)}$ and $\sigma =
\frac{np - rn + pr}{p}$. The second step combined with (\ref{3int})
provide

\[
\int_{\tilde{D}_\alpha} u_\alpha^r dv_g \rightarrow 0
\]

\n as $\alpha \rightarrow \infty$, so that $m_\alpha \rightarrow
\infty$. We now derive a contradiction from (\ref{3dcm}). From
(\ref{3lu}), (\ref{3dqp}) and (\ref{dt}), we find

\begin{equation}\label{3dpt}
m_\alpha^\sigma \int_{\tilde{D}_\alpha} u_\alpha^r dv_g \leq
m_\alpha^\sigma ||u_\alpha||^r_{L^\infty(\tilde{D}_\alpha)}
(A_\alpha^{\frac{1}{p}} u_\alpha(y_\alpha)^{\frac{p-r}{p}})^n \leq c
m_\alpha^\sigma u_\alpha(y_\alpha)^r (A_\alpha^{\frac{1}{p}}
u_\alpha(y_\alpha)^{\frac{p-r}{p}})^n \leq c\ .
\end{equation}

\n Consider the function $\eta_\alpha(x) = \eta( A_\alpha^{-1/p}
d_g(x,y_\alpha) u_\alpha(y_\alpha)^{\frac{r-p}{p}})$, where $\eta
\in C_0^1(\R)$ is a cutoff function such that $\eta = 1$ on
$[0,\frac{1}{2}]$, $\eta = 0 $ on $[1,\infty)$ and $0 \leq \eta \leq
1$. Taking $u_\alpha \eta_\alpha^p$ as a test function in
(\ref{3ep}), one has

\[
A_\alpha \int_M |\nabla_g u_\alpha|^p \eta_\alpha^p dv_g + p
A_\alpha \int_M |\nabla_g u_\alpha|^{p - 2} u_\alpha \eta_\alpha^{p
- 1} \nabla_g u_\alpha \cdot \nabla_g \eta_\alpha dv_g + \alpha
A_\alpha \int_M u_\alpha^p \eta_\alpha^p dv_g
\]
\[
+ \frac{1 - \theta}{\theta} B_\alpha \int_M u_\alpha^q \eta_\alpha^p
dv_g = \mu_\alpha \int_M u_\alpha^r \eta_\alpha^p dv_g\ .
\]

\n Clearly,

\[
\left| \int_M |\nabla_g u_\alpha|^{p - 2} u_\alpha \eta_\alpha^{p -
1} \nabla_g u_\alpha \cdot \nabla_g \eta_\alpha\; dv_g \right| \leq
\varepsilon \int_M |\nabla_g u_\alpha|^p \eta_\alpha^p dv_g +
c_\varepsilon \int_{\tilde{D}_\alpha} |\nabla_g \eta_\alpha|^p
u_\alpha^p dv_g\ .
\]

\n From (\ref{3lu}), (\ref{3dqp}) and (\ref{dt}), it follows that

\begin{equation}\label{3dm}
A_\alpha \int_M|\nabla_g \eta_\alpha|^p u_\alpha^p dv_g \leq
A_\alpha (A_\alpha^{-1/p} u_\alpha(y_\alpha)^{\frac{r-p}{p}})^p
\int_{\tilde{D}_\alpha} u_\alpha^p dv_g \leq c u_\alpha(y_\alpha)^r
(A_\alpha^{\frac{1}{p}} u_\alpha(y_\alpha)^{\frac{p-r}{p}})^n \leq c
m_\alpha^{- \sigma}\ .
\end{equation}

\n Consequently, these inequalities, (\ref{3dpt}) and (b) of the
first step imply

\begin{equation}\label{31d}
A_\alpha \int_M |\nabla_g u_\alpha|^p \eta_\alpha^p dv_g + c \alpha
A_\alpha \int_M u_\alpha^p \eta_\alpha^p dv_g + c B_\alpha \int_M
u_\alpha^q \eta_\alpha^p dv_g \leq c m_\alpha^{-\sigma}\ .
\end{equation}

\n On the other hand, \ref{AB1} provides

\[
\left( \int_{\tilde{D}_\alpha} u_\alpha^r dv_g \right)^{\frac{p}{r
\theta}} \leq \left( \int_M (u_\alpha \eta_\alpha^p)^r dv_g
\right)^{\frac{p}{r \theta}} \leq c \left(\int_M |\nabla_g
u_\alpha|^p \eta_\alpha^{p^2} dv_g\right) \left( \int_M (u_\alpha
\eta_\alpha^p)^q dv_g \right)^{\frac{p(1 - \theta)}{\theta q}}
\]

\begin{equation}\label{3ddb}
+ c \left(\int_M |\nabla_g \eta_\alpha|^p u_\alpha^p dv_g\right)
\left(\int_M (u_\alpha \eta_\alpha^p)^q dv_g \right)^{\frac{p(1 -
\theta)}{\theta q}} + c \left(\int_M (u_\alpha \eta_\alpha^p)^p
dv_g\right) \left( \int_M (u_\alpha \eta_\alpha^p)^q dv_g
\right)^{\frac{p(1 - \theta)}{\theta q}}\ .
\end{equation}

\n Using (\ref{3dm}) and (\ref{31d}), we may then estimate each term of the
right-hand side of (\ref{3ddb}). By (a), (c) and (d) of the
first step,

\[
\left(\int_M |\nabla_g u_\alpha|^p \eta_\alpha^{p^2} dv_g\right)
\left( \int_M (u_\alpha \eta_\alpha^p)^q dv_g \right)^{\frac{p(1 -
\theta)}{\theta q}}
\]

\[
\leq  \left(\frac{A_\alpha \int_M |\nabla_g u_\alpha|^p
\eta_\alpha^p dv_g}{A_\alpha B_\alpha^{\frac{p(1 - \theta)}{\theta
q}}}\right) \left( B_\alpha \int_M u_\alpha^q \eta_\alpha^p dv_g
\right)^{\frac{p(1 - \theta)}{\theta q}} \leq c m_\alpha^{-\sigma(1
+ \frac{p( 1 - \theta)}{\theta q})},
\]

\[
\left(\int_M |\nabla_g \eta_\alpha|^p u_\alpha^p dv_g\right)
\left(\int_M (u_\alpha \eta_\alpha^p)^q dv_g \right)^{\frac{p(1 -
\theta)}{\theta q}}
\]

\[
\leq \frac{A_\alpha \int_M |\nabla_g \eta_\alpha|^p u_\alpha^p
dv_g}{A_\alpha B_\alpha^{\frac{p( 1 - \theta)}{\theta q}}} \left(
B_\alpha \int_M u_\alpha^q \eta_\alpha^p dv_g \right)^{\frac{p(1 -
\theta)}{\theta q}} \leq c m_\alpha^{-\sigma(1 + \frac{p(1 -
\theta)}{\theta q})}
\]

\n and

\[
\left(\int_M (u_\alpha \eta_\alpha^p)^p dv_g\right) \left( \int_M
(u_\alpha \eta_\alpha^p)^q dv_g \right)^{\frac{p(1 - \theta)}{\theta
q}}
\]

\[
\leq \frac{A_\alpha \int_M u_\alpha^p \eta_\alpha^p dv_g}{A_\alpha
B_\alpha^{\frac{p(1 - \theta)}{\theta q}}} \left(B_\alpha \int_M
\eta_\alpha^p u_\alpha^q dv_g \right)^{\frac{p(1 - \theta)}{\theta
q}} \leq c m_\alpha^{-\sigma(1 + \frac{p(1 - \theta)}{\theta q})}\ .
\]

\n Replacing these estimates in (\ref{3ddb}), one has

\[
\left(\int_{\tilde{D}_\alpha} u_\alpha^r dv_g \right)^{\frac{p}{r
\theta}} \leq c m_\alpha^{-\sigma(1 + \frac{p(1 - \theta)}{\theta
q})},
\]

\n so that

\[
m_\alpha^\sigma \int_{\tilde{D}_\alpha} u_\alpha^r dv_g \leq c
m_\alpha^{\sigma( -\frac{r \theta}{p}(1 + \frac{p(1 -
\theta)}{\theta q}) + 1)}\ .
\]

\n From the range of $p$, $q$ and $r$, it follows easily that

\[
- \frac{r \theta}{p} \left( 1 + \frac{p(1 - \theta)}{\theta
q}\right) + 1 < 0\ .
\]

\n So,

\[
m_\alpha^\sigma \int_{\tilde{D}_\alpha} u_\alpha^r dv_g \rightarrow
0
\]

\n since $m_\alpha \rightarrow \infty$. But this contradicts
(\ref{3dcm}), and so the part (a) is proved.

In order to prove the part (b), we again argue by contradiction.
Since part of the arguments are similar to those ones already used, we
will omit some details. Set $w_\alpha(x) =
d_g(x,x_\alpha)^{\frac{n}{p}} u_\alpha(x)$. Let $y_\alpha \in
\overline{B(x_\alpha,1)}$ be a maximum point of $w_\alpha$. Then,
$w_\alpha(y_\alpha) = ||w_\alpha||_{L^\infty(B(x_\alpha,1))}
\rightarrow \infty$. Arguing as in (a), one has

\begin{equation}\label{n2d7}
B(y_\alpha,u_\alpha(y_\alpha)^{-\frac{p}{n}}) \cap B(x_\alpha,
A_\alpha^{\frac{1}{p}} w_\alpha(y_\alpha)^\nu) = \emptyset
\end{equation}

\n for some $\nu > 0$ and $\alpha$ large. From $w_\alpha(y_\alpha)
\rightarrow \infty$, it follows easily that $a_\alpha =
||u_\alpha||_{L^\infty(M)}^{- \frac{n + p^2}{pn}}
u_\alpha(y_\alpha)^{\frac{1}{p}} \rightarrow 0$. For $x \in
B(0,1)$, we set

\[
\begin{array}{l}
h_\alpha(x) = g(\exp_{y_\alpha}(a_\alpha x)),\vspace{0,3cm}\\
\psi_\alpha(x) = u_\alpha(y_\alpha)^{-1} u_\alpha(\exp_{y_\alpha}(
a_\alpha x))\ .
\end{array}
\]

\n Clearly, $\psi_\alpha$ satisfies

\[
\Delta_{p, h_\alpha} \psi_\alpha + \alpha a_\alpha^p \psi_\alpha^{p
- 1} + \frac{1 - \theta}{\theta} \frac{a_\alpha^p B_\alpha
u_\alpha(y_\alpha)^{q - p}}{A_\alpha} \psi_\alpha^{q - 1} =
\mu_\alpha \frac{a_\alpha^p}{A_\alpha} \psi_\alpha^{p - 1}.
\]

\n We now verify that each coefficient of this equation is bounded.
By (\ref{3dqp}) and (c) of the first step,

\[
\alpha a_\alpha^p = \alpha ||u_\alpha||_{L^\infty(M)}^{- \frac{n +
p^2}{n}} u_\alpha(y_\alpha) \leq c \alpha A_\alpha \leq c\ .
\]

\n By (d) of the first step, one has

\[
B_\alpha ||u_\alpha||_{L^\infty(M)}^{q - p} \leq c,
\]

\n so that, by (\ref{3dqp}),

\[
\frac{a_\alpha^p B_\alpha u_\alpha(y_\alpha)^{q - p}}{A_\alpha} \leq
c \frac{||u_\alpha||_{L^\infty(M)}^{-\frac{n + p^2}{n} + p -q}
u_\alpha(y_\alpha)^{q - p + 1}}{A_\alpha} \leq c
\frac{||u_\alpha||_{L^\infty(M)}^{- \frac{p^2}{n}}}{A_\alpha} \leq c
\]

\n since $q - p + 1 \geq 0$, and

\[
\frac{a_\alpha^p}{A_\alpha} = \frac{||u_\alpha||_{L^\infty(M)}^{-
\frac{n + p^2}{n}} u_\alpha(y_\alpha)}{A_\alpha} \leq c\ .
\]

\n Arguing as in (\ref{3lu}), one easily checks

\begin{equation}\label{dps}
||u_\alpha||_{L^\infty(B(y_\alpha,u_\alpha(y_\alpha)^{-\frac{p}{n}}))}
\leq c u_\alpha(y_\alpha)\ .
\end{equation}

\n Hence, we have

\[
||\psi_\alpha||_{L^{\infty}(B(0,1))} \leq c \ .
\]

\n By Tolksdorf \cite{To}, we conclude that $\psi_\alpha$
converges to $\psi$ in $C_{loc}^1(B(0,1))$. Noting that $\psi \not
\equiv 0$, by (\ref{3dqp}),

\[
0 < c \leq \int_{B(0,1)} \psi_\alpha^q dv_{h_\alpha} =
u_\alpha(y_\alpha)^{-q} a_\alpha^{-n} A_\alpha^{\frac{\theta q}{p(1
- \theta)}} \frac{\int_{B(y_\alpha,a_\alpha)} u_\alpha^q
dv_g}{\int_M u_\alpha^q dv_g}
\]

\begin{equation}\label{3ctba}
\leq c
\left(\frac{||u_\alpha||_{L^\infty(M)}}{u_\alpha(y_\alpha)}\right)^{
\frac{n + pq}{p}} \frac{\int_{B(y_\alpha,a_\alpha)} u_\alpha^q
dv_g}{\int_M u_\alpha^q dv_g}\ .
\end{equation}

\n Therefore, Lemma \ref{3lmi} and (\ref{n2d7}) imply

\[
m_\alpha = \frac{||u_\alpha||_{L^\infty(M)}}{u_\alpha(y_\alpha)}
\rightarrow \infty\ .
\]

\n We now show by induction that exists $\gamma_k \rightarrow
\infty$ such that, for each $k$,

\begin{equation}\label{n2d8}
m_\alpha^{\gamma_k} \int_{B(y_\alpha,2^{-k}
u_\alpha(y_\alpha)^{-\frac{p}{n}})} u_\alpha^p dv_g \rightarrow 0 \
.
\end{equation}

\n By (\ref{3dqp}) and (\ref{dps}),

\[
\int_{B(y_\alpha,u_\alpha(y_\alpha)^{-\frac{p}{n}})} u_\alpha^p dv_g
\leq c m_\alpha^{p - q}
\frac{\int_{B(y_\alpha,u_\alpha(y_\alpha)^{-\frac{p}{n}})}
u_\alpha^q dv_g}{\int_M u_\alpha^q dv_g}\ .
\]

\n Then, applying Lemma \ref{3lmi} and (\ref{n2d7}), one obtains

\[
m_\alpha^{p - q}
\int_{B(y_\alpha,u_\alpha(y_\alpha)^{-\frac{p}{n}})} u_\alpha^p dv_g
\leq c \frac{\int_{B(y_\alpha,u_\alpha(y_\alpha)^{-\frac{p}{n}})}
u_\alpha^q dv_g}{\int_M u_\alpha^q dv_g} \rightarrow 0\ .
\]

\n Let $\gamma_0 = p - q$. Supposing that (\ref{n2d8}) holds for
some $\gamma_k \geq \gamma_0$, we will then show that (\ref{n2d8}) also holds for
$\gamma_{k + 1} = \gamma_k(1 + \nu)$, where $\nu > 0$ is a fixed
positive number to be determined later. Consider the function
$\eta_{\alpha,k}(x) = \eta(2^k u_\alpha(y_\alpha)^{\frac{p}{n}}
d_g(y_\alpha,x))$, where $\eta \in C_0^1(\R)$ satisfies $\eta = 1$
on $[0,\frac{1}{2}]$, $\eta = 0 $ on $[1,\infty)$ and $0 \leq \eta
\leq 1$. Taking $u_\alpha \eta_{\alpha,k}^p$ as a test function in
(\ref{3ep}) and using the first step, one finds

\[
A_\alpha \int_M |\nabla_g u_\alpha|^p \eta_{\alpha,k}^p dv_g +
\alpha A_\alpha \int_M u_\alpha^p \eta_{\alpha,k}^p dv_g + B_\alpha
\int_M u_\alpha^q \eta_{\alpha,k}^p dv_g
\]

\[
\leq c \int_M u_\alpha^p \eta_{\alpha,k}^p dv_g + c A_\alpha \int_M
|\nabla_g \eta_{\alpha,k}|^p u_\alpha^p dv_g\ .
\]

\n By (\ref{3dqp}) and (\ref{n2d8}), one has

\[
\int_M u_\alpha^p \eta_{\alpha,k}^p dv_g \leq c
m_\alpha^{-\gamma_k}
\]

\n and

\begin{equation}\label{n2d9}
A_\alpha \int_M |\nabla_g \eta_{\alpha,k}|^p u_\alpha^p dv_g \leq c
m_\alpha^{-\gamma_k}
\end{equation}

\n since $|\nabla_g \eta_{\alpha,k}| \leq c u_\alpha(y_\alpha)^{\frac{p}{n}}$. Then,

\begin{equation}\label{n2d10}
A_\alpha \int_M |\nabla_g u_\alpha|^p \eta_{\alpha,k}^p dv_g +
\alpha A_\alpha \int_M u_\alpha^p \eta_{\alpha,k}^p dv_g + B_\alpha
\int_M u_\alpha^q \eta_{\alpha,k}^p dv_g \leq c
m_\alpha^{-\gamma_k}\ .
\end{equation}

\n On the other hand, using \ref{AB1}, one has

\[
\int_{B(y_\alpha, 2^{-(k + 1)} u_\alpha(y_\alpha)^{- \frac{p}{n}})}
u_\alpha^p dv_g \leq \int_M(u_\alpha \eta_{\alpha,k}^p)^p dv_g
\]

\[
\leq c \left( \left( \int_M |\nabla_g u_\alpha|^p \eta_{\alpha,k}^p
dv_g + \int_M |\nabla_g \eta_{\alpha,k}|^p u_\alpha^p dv_g + \int_M
u_\alpha^p \eta_{\alpha,k}^p dv_g \right)\left(\int_M u_\alpha^q
\eta_{\alpha,k}^p dv_g \right)^{\frac{p(1 - \theta)}{\theta
q}}\right)^\theta\ .
\]

\n Therefore, (\ref{n2d9}), (\ref{n2d10}) and (d) of the first step
imply

\[
\int_{B(y_\alpha,2^{-(k+1)} u_\alpha(y_\alpha)^{-\frac{p}{n}})}
u_\alpha^p dv_g \leq c m_\alpha^{-\gamma_k (1 + \frac{p(1 -
\theta)}{\theta q})\theta}\ .
\]

\n From the definition of $\theta$, it follows that this estimate may be rewritten as

\[
\int_{B(y_\alpha,2^{-(k+1)} u_\alpha(y_\alpha)^{-\frac{p}{n}})}
u_\alpha^p dv_g \leq c m_\alpha^{-\gamma_{k+1}},
\]

\n where

\[
\gamma_{k+1} = \gamma_k\left( 1 + \frac{p(p - q)}{n(p - q) + pq}
\right) \ .
\]

\n So, the proof follows by induction and, moreover,
$\gamma_k \rightarrow \infty$. Fixing $k$ large such that
$\gamma_k \geq \frac{n + p^2}{p}$, we arrive at the following
contradiction

\[
0 < c \leq \int_{B(0,2^{-k})} \psi_\alpha^p dv_{h_\alpha} \leq
u_\alpha(y_\alpha)^{-\frac{n + p^2}{p}}
||u_\alpha||_{L^\infty(M)}^{\frac{n + p^2}{p}} \int_{B(y_\alpha,
2^{-k} u_\alpha (y_\alpha)^{-\frac{p}{n}})} u_\alpha^p dv_g \leq c
m_\alpha^{\frac{n + p^2}{p} - \gamma_k} \rightarrow 0\ .
\]

\n This ends the proof of the third step. \bl \\

\n {\bf Final conclusions:} In this last part, we will estimate several
integrals by using the third step and will derive a contradiction at the
end. Let $\eta \in C^1_0(\R)$ be a function such that $\eta = 1$ on
$[0,1)$, $\eta = 0$ on $[2, \infty)$ and $0 \leq \eta \leq 1$.
Define $\eta_\alpha(x) = \eta(d_g(x,x_\alpha))$. The optimal
Euclidean Gagliardo-Nirenberg inequality provides

\[
\left( \int_{B(x_\alpha,1)} u_\alpha^r dx \right)^{\frac{p}{r
\theta}} \leq \left( \int_{B(x_\alpha, 2)} (u_\alpha \eta_\alpha)^r
dx \right)^{\frac{p}{r \theta}} \leq A(p,q,r) \left(
\int_{B(x_\alpha, 2)} |\nabla (u_\alpha \eta_\alpha)|^p dx\right)
\left(\int_{B(x_\alpha, 2)} (u_\alpha \eta_\alpha)^q dx
\right)^{\frac{p(1 - \theta)}{\theta q}}\ .
\]

\n Using normal coordinates and the Cartan expansion

\[
(1 - c d_g(x,x_\alpha)^2) dv_g \leq dx \leq (1 + c
d_g(x,x_\alpha)^2) dv_g,
\]
\n we obtain

\[
\left( \int_{B(x_\alpha, 1)} u_\alpha^r dx
\right)^{\frac{p}{r\theta}} \leq A(p,q,r) \left( \int_{B(x_\alpha,
2)} |\nabla_g (u_\alpha \eta_\alpha)|^p dv_g + c \int_{B(x_\alpha,
2)} |\nabla_g (u_\alpha\eta_\alpha)|^p d_g(x,x_\alpha)^2dv_g \right)
\]

\[ \times \left(A_\alpha + c \left( \int_M u_\alpha^q dv_g
\right)^{\frac{p(1 - \theta)}{\theta q} - 1} \left( \int_M (u_\alpha
\eta_\alpha)^q d_g(x,x_\alpha)^2 dv_g\right) \right)\ .
\]

\n Applying then the inequalities

\[
|\nabla_g (u_\alpha \eta_\alpha)|^p \leq |\nabla_g u_\alpha|^p + c
|\eta_\alpha \nabla_g u_\alpha|^{p - 1} |u_\alpha \nabla_g
\eta_\alpha| + c |u_\alpha \nabla_g \eta_\alpha|^p
\]

\n and

\[
A(p,q,r) \left( A_\alpha \int_M |\nabla_g u_\alpha|^p dv_g\right) +
A(p,q,r) \left( \alpha A_\alpha \int_M u_\alpha^p dv_g \right) < 1,
\]

\n one arrives at

\[
A(p,q,r) \alpha A_\alpha \int_M u_\alpha^p dv_g \leq 1 - \left(
\int_{B(x_\alpha, 1)} u_\alpha^r dx \right)^{\frac{p}{r \theta}} + c
A_\alpha \int_M u_\alpha^p dv_g + c A_\alpha \int_{B(x_\alpha, 2)}
\eta_\alpha^p |\nabla_g u_\alpha|^p d_g(x,x_\alpha)^2 dv_g
\]

\[
 + c \left( \int_{B(x_\alpha,2)} |\nabla_g u_\alpha|^p dv_g\right)
\left(\int_M u_\alpha^q dv_g\right)^{\frac{p(1 - \theta)}{\theta q}
-1} \left( \int_M (u_\alpha \eta_\alpha)^q d_g(x,x_\alpha)^2 dv_g
\right)
\]

\begin{equation}\label{3du}
 + c A_\alpha \int_{B(x_\alpha, 2)} \eta_\alpha^{p - 1} |\nabla_g
u_\alpha|^{p -1} |\nabla_g \eta_\alpha| u_\alpha dv_g\ .
\end{equation}

\n We now estimate each term of the right-hand side of (\ref{3du}).
By (\ref{fim31}), we may apply the second step with $c_\alpha = 1$ for all $\alpha$, so
that $\int_{B(x_\alpha,\delta)} u_\alpha^r dv_g \rightarrow 1$. In particular, there exists $c>0$ such that

\[
c^{-1} \leq \int_{B(x_\alpha,1)} u_\alpha^r dx \leq c\ .
\]

\n So, the Cartan expansion and $||u_\alpha||_{L^r(M)} = 1$ lead us
to

\begin{equation}\label{est1}
1 - \left(\int_{B(x_\alpha, 1)} u_\alpha^r dx \right)^{\frac{p}{r
\theta}} \leq c \left(1 - \int_{B(x_\alpha, 1)} u_\alpha^r dx
\right) \leq c \int_{M \backslash B(x_\alpha, 1)} u_\alpha^r dv_g +
c \int_{B(x_\alpha, 1)} u_\alpha^r d_g(x,x_\alpha)^2 dv_g\ .
\end{equation}

\n Let $\zeta \in C^1(\R)$ be a function such that $0 \leq \zeta
\leq 1$, $\zeta = 0$ on $[0, \frac{1}{2}]$, $\zeta = 1$ on $[1,
\infty)$. Define $\zeta_\alpha = \zeta(d_g(x,x_\alpha))$. Taking
$u_\alpha \zeta_\alpha^p$ as a test function in (\ref{3ep}), one finds

\[
A_\alpha \int_M \zeta_\alpha^p |\nabla_g u_\alpha|^p dv_g \leq c
\int_M \zeta_\alpha^p u_\alpha^r dv_g + c  A_\alpha \int_M
\zeta_\alpha^{p - 1} |\nabla_g u_\alpha|^{p-1} |\nabla_g
\zeta_\alpha| u_\alpha dv_g,
\]

\n so that

\begin{equation}\label{3dfrq}
A_\alpha \int_M \zeta_\alpha^p |\nabla_g u_\alpha|^p dv_g \leq c
A_\alpha \int_M u_\alpha^p dv_g + c \int_{M \backslash
B(x_\alpha,\frac{1}{2})} u_\alpha^r dv_g\ .
\end{equation}

\n This implies

\begin{equation}\label{3dfq1}
A_\alpha \int_M \zeta_\alpha^p |\nabla_g u_\alpha|^{p - 1} u_\alpha
dv_g \leq c A_\alpha \int_M u_\alpha^p dv_g + c \int_{M \setminus
B(x_\alpha, \frac{1}{2})} u_\alpha^r dv_g\ .
\end{equation}

\n Since $p \leq 2$, it follows that

\begin{equation}\label{3dfq2}
\int_M u_\alpha \eta_\alpha^p |\nabla_g u_\alpha|^{p - 1}
d_g(x,x_\alpha) dv_g \leq \varepsilon \int_M \eta_\alpha^p |\nabla_g
u_\alpha|^p d_g(x,x_\alpha)^2 dv_g + c_{\varepsilon} \int_M
u_\alpha^p dv_g\ .
\end{equation}

\n Taking $u_\alpha d_g^2 \eta_\alpha^p$ as a test function in
(\ref{3ep}), one has

\[
A_\alpha \int_M \eta_\alpha^p |\nabla_g u_\alpha|^p
d_g(x,x_\alpha)^2 dv_g \leq c \int_{B(x_\alpha, 2)} u_\alpha^r
d_g(x,x_\alpha)^2 dv_g + c A_\alpha \int_M u_\alpha \eta_\alpha^p
|\nabla_g u_\alpha|^{p - 1} d_g(x,x_\alpha) dv_g
\]

\[
+ c A_\alpha\int_M \zeta_\alpha^p |\nabla_g u_\alpha|^{p - 1}
u_\alpha dv_g\ .
\]

\n So, from (\ref{3dfq1}) and (\ref{3dfq2}), we obtain

\begin{equation} \label{est2}
A_\alpha \int_M \eta_\alpha^p |\nabla_g u_\alpha|^p
d_g(x,x_\alpha)^2 dv_g \leq c \int_{B(x_\alpha, 2)} u_\alpha^r
d_g(x,x_\alpha)^2 dv_g + c \int_{M \backslash B(x_\alpha,
\frac{1}{2})} u_\alpha^r dv_g + c A_\alpha \int_M u_\alpha^p dv_g\ .
\end{equation}

\n It follows from (\ref{3dfrq}) that

\begin{equation}\label{est3}
A_\alpha \int_{B(x_\alpha, 2)} \eta_\alpha^{p - 1} |\nabla_g
u_\alpha|^{p -1} |\nabla_g \eta_\alpha| u_\alpha dv_g \leq  c
A_\alpha \int_M u_\alpha^p dv_g + c \int_{M \backslash B(x_\alpha,
1)} u_\alpha^r dv_g\ .
\end{equation}

\n Finally, taking $u_\alpha d_g^2 \eta_\alpha^p$ as a test function
in (\ref{3ep}), one has

\[
B_\alpha \int_M u_\alpha^q \eta_\alpha^p d_g(x,x_\alpha)^2 dv_g \leq
c \int_{B(x_\alpha, 2)}u^r_\alpha d_g(x,x_\alpha)^2 dv_g + A_\alpha
\int_{B(x_\alpha, 2)} \eta_\alpha^{p - 1} |\nabla_g u_\alpha|^{p -1}
|\nabla_g \eta_\alpha| u_\alpha dv_g
\]

\begin{equation}\label{3dp4}
+ \int_M u_\alpha \eta_\alpha^p |\nabla_g u_\alpha|^{p - 1}
d_g(x,x_\alpha) dv_g\ .
\end{equation}

\n Using (\ref{est3}), (\ref{3dfq2}) and (\ref{est2}), one arrives
at

\[
B_\alpha \int_M u_\alpha^q \eta_\alpha^p d_g(x,x_\alpha)^2 dv_g \leq
c A_\alpha \int_M u_\alpha^p dv_g + c \int_{B(x_\alpha,
2)}u^r_\alpha d_g(x,x_\alpha)^2 dv_g + c \int_{M \setminus
B(x_\alpha,\frac{1}{2})} u_\alpha^r dv_g,
\]

\n or equivalently,

\[
\left( \int_{B(x_\alpha,2)} |\nabla_g u_\alpha|^p dv_g\right)
\left(\int_M u_\alpha^q dv_g\right)^{\frac{p(1 - \theta)}{\theta q}
-1} \left( \int_M (u_\alpha \eta_\alpha)^q d_g(x,x_\alpha)^2 dv_g
\right) \leq c A_\alpha \int_M u_\alpha^p dv_g
\]

\begin{equation} \label{est4}
+ c \int_{B(x_\alpha, 2)}u^r_\alpha d_g(x,x_\alpha)^2 dv_g + c
\int_{M \setminus B(x_\alpha,\frac{1}{2})} u_\alpha^r dv_g\ .
\end{equation}

\n Replacing (\ref{est1}), (\ref{est2}), (\ref{est4}) and
(\ref{est3}) in (\ref{3du}), we obtain

\begin{equation}\label{ud}
\alpha A_\alpha \int_M u_\alpha^p dv_g \leq c A_\alpha \int_M
u_\alpha^p dv_g + c \int_{M \backslash B(x_\alpha, \frac{1}{2})}
u_\alpha^r dv_g + c \int_{B(x_\alpha, 2)} u_\alpha^r
d_g(x,x_\alpha)^2 dv_g\ .
\end{equation}

\n We now analyze separately the cases $p < r$ and $p = r$. For $p
< r$, we use the part (a) of the third step and obtain

\[
\int_{M\backslash B(x_\alpha, \frac{1}{2})} u_\alpha^r dv_g =
\int_{M \backslash B(x_\alpha, \frac{1}{2})} d_g(x,x_\alpha)^{-p}
d_g(x,x_\alpha)^p u_\alpha^{r - p} u_\alpha^p dv_g \leq c A_\alpha
\int_M u_\alpha^p dv_g
\]

\n and

\[
\int_{B(x_\alpha,2)} u_\alpha^r d_g(x,x_\alpha)^2 dv_g =
\int_{B(x_\alpha, 2)} d_g(x,x_\alpha)^{2 - p} d_g(x,x_\alpha)^p
u_\alpha^{r - p} u_\alpha^p dv_g \leq c A_\alpha \int_M u_\alpha^p
dv_g
\]

\n since $p \leq 2$. Replacing these inequalities in (\ref{ud}), one
finds $\alpha \leq c$. This contradiction ends the proof of
Theorem \ref{tgno1} in the case $p < r$. Finally, consider the case
$p = r$. From the part (b) of the third step, one has

\[
\int_{M \setminus B(x_\alpha,\frac{1}{2})} u_\alpha^p dv_g \leq
\int_{M \setminus B(x_\alpha,\frac{1}{2})} d_g(x,x_\alpha)^{-p}
d_g(x,x_\alpha)^p (u_\alpha^{\frac{p}{n}})^p u_\alpha^{\frac{np -
p^2}{n}} dv_g \leq \int_M u_\alpha^{\frac{np - p^2}{n}} dv_g,
\]

\n and

\begin{equation} \label{est-d}
\int_{B(x_\alpha,2)} u_\alpha^p d_g(x,x_\alpha)^2 dv_g =
\int_{B(x_\alpha, 2)} d_g(x,x_\alpha)^{2 - p} d_g(x,x_\alpha)^p
(u_\alpha^{\frac{p}{n}})^p u_\alpha^{\frac{np - p^2}{n}} dv_g \leq c
\int_M u_\alpha^{\frac{np - p^2}{n}} dv_g\ .
\end{equation}

\n If $q \leq \frac{np - p^2}{n}$, by interpolation inequality, one
concludes that

\[
\int_M u_\alpha^{\frac{np - p^2}{n}} dv_g \leq \left( \int_M
u_\alpha^q dv_g \right)^{\frac{p^2}{n(p-q)}} = A_\alpha \ .
\]

\n In this case, from (\ref{ud}), we again obtain the contradiction
$\alpha \leq c$. For $q \geq \frac{np - p^2}{n}$, we have

\[
\int_M u_\alpha^{\frac{np - p^2}{n}} dv_g \leq c \left( \int_M
u_\alpha^q dv_g \right)^{\frac{np - p^2}{qn}} = c A_\alpha^{\frac{(n
- p)(p - q)}{pq}}\ .
\]

\n Set $\delta_ 0 = \frac{(n - p)(p - q)}{pq}$. Let $\zeta \in
C^1(\R)$ be a function such that $0 \leq \zeta \leq 1$, $\zeta = 0$
on $[0, \frac{1}{2}]$, $\zeta = 1$ on $[1, \infty)$. Define
$\zeta_\alpha = \zeta(d_g(x,x_\alpha))$. Taking $u_\alpha
\zeta_\alpha^p$ as a function test in (\ref{3ep}) and arguing in a
similar manner to the proof of (\ref{n2d8}), one arrives at a constant $c > 0$ such that

\[
\int_{M \setminus B(x_\alpha,\frac{1}{2})} u_\alpha^p dv_g \leq c
A_\alpha^{\delta_k},
\]

\n where $\delta_k$ satisfies the recurrence relation, for any $k
\geq 0$,

\[
\delta_{k+1} = \delta_k (1 + \frac{p(p - q)}{q(p - n) + np})\ .
\]

\n In particular, $\delta_k \geq 1$ for $k$ large. We now show that

\[
\int_{B(x_\alpha,2)} u_\alpha^p d_g(x,x_\alpha)^2 dv_g \leq c
A_\alpha^{\delta_k}
\]

\n for all $k \geq 0$. From (\ref{est-d}), we have

\[
\int_{B(x_\alpha,2)} u_\alpha^p d_g(x,x_\alpha)^2 dv_g \leq c
A_\alpha^{\delta_0}\ .
\]

\n Let $\eta_\alpha(x) = \eta(d_g(x,x_\alpha)$, where $\eta \in
C_0^1(\R)$ satisfies $\eta = 1$ on $[0,\frac{1}{2}]$, $\eta = 0$ on
$[1,\infty)$ and $0 \leq \eta \leq 1$. Taking $u_\alpha
\eta_\alpha^p d_g^p$ as a test function in (\ref{3ep}), we obtain

\[
A_\alpha \int_M |\nabla_g u_\alpha|^p \eta_{\alpha}^p
d_g(x,x_\alpha)^p dv_g + \alpha A_\alpha \int_M u_\alpha^p
\eta_{\alpha}^p d_g(x,x_\alpha)^p dv_g + B_\alpha \int_M u_\alpha^q
\eta_{\alpha}^p d_g(x,x_\alpha)^p dv_g
\]

\[
\leq c A_\alpha^{\delta_0} + A_\alpha \int_M |\nabla_g
u_\alpha|^{p-1} \eta_{\alpha}^{p-1} u_\alpha |\nabla_g
\eta_{\alpha}| d_g(x,x_\alpha)^p dv_g + A_\alpha \int_M |\nabla_g
u_\alpha|^{p-1} \eta_{\alpha}^p u_\alpha d_g(x,x_\alpha)^{p - 1}
dv_g,
\]

\n so that, applying H\"{o}lder and Young inequalities,

\[
A_\alpha \int_M |\nabla_g u_\alpha|^p \eta_{\alpha}^p
d_g(x,x_\alpha)^p dv_g + \alpha A_\alpha \int_M u_\alpha^p
\eta_{\alpha}^p d_g(x,x_\alpha)^p dv_g + B_\alpha \int_M u_\alpha^q
\eta_{\alpha}^p d_g(x,x_\alpha)^p dv_g
\]

\[
\leq c A_\alpha^{\delta_0} + c A_\alpha \int_{M} u_\alpha^p dv_g
\leq c A_\alpha^{\delta_0}\ .
\]

\n On the other hand, \ref{AB1} and the hypothesis $\frac{p^2}{2} \leq q$ provide

\[
\left(\int_{B(x_\alpha,\frac12)} u_\alpha^p d_g(x,x_\alpha)^2
dv_g\right)^{\frac{1}{\theta}} \leq \left(\int_M (u_\alpha
\eta_\alpha^p d_g(x,x_\alpha)^{\frac{2}{p}})^p
dv_g\right)^{\frac{1}{\theta}}
\]

\[
\leq c \left( \int_M |\nabla_g u_\alpha|^p \eta_{\alpha}^p
d_g(x,x_\alpha)^2 dv_g + \int_M |\nabla_g u_\alpha|^{p-1}
\eta_{\alpha}^{p-1} d_g(x,x_\alpha)^{\frac{2(p-1)}{p}} u_\alpha
|\nabla_g (\eta_{\alpha}^p d_g(x,x_\alpha)^{\frac2p})| dv_g \right.
\]
\[
\left. + \int_M u_\alpha^p |\nabla_g (\eta_{\alpha}^p
d_g(x,x_\alpha)^{\frac2p})|^p dv_g + \int_M u_\alpha^p
d_g(x,x_\alpha)^2 dv_g \right) \left( \int_M u_\alpha^q
\eta_\alpha^p d_g(x,x_\alpha)^{\frac{2q}{p}}  dv_g
\right)^{\frac{p(1-\theta)}{\theta q}}
\]
\[
\leq c \left( \int_M |\nabla_g u_\alpha|^p \eta_{\alpha}^p
d_g(x,x_\alpha)^2 dv_g + \int_{M \setminus B(x_\alpha, \frac12)}
u_\alpha^p dv_g + 1 \right) \left( \int_M u_\alpha^q \eta_\alpha^p
d_g(x,x_\alpha)^p dv_g \right)^{\frac{p(1-\theta)}{\theta q}} \leq c
(A_\alpha^{\delta_0(1 + \nu)})^{\frac{1}{\theta}}\ .
\]

\n Hence,

\[
\int_{B(x_\alpha,\frac12)} u_\alpha^p d_g(x,x_\alpha)^2 dv_g \leq c
A_\alpha^{\delta_1}\ .
\]

\n Arguing by induction, it easily follows that

\[
\int_{B(x_\alpha,1)} u_\alpha^p d_g(x,x_\alpha)^2 dv_g \leq c
A_\alpha^{\delta_k}
\]

\n for all $k \geq 0$. Combining this estimate with (\ref{ud}) we
again arrive at the contradiction $\alpha
\leq c$. This ends the proof of Theorem \ref{tgno1}. \bl \\

{\bf Acknowledgements.} The authors thank the referee for his valuable comments. The first author was partially supported by Brazilian Agency CNPq.

\end{document}